\documentclass[graybox]{svmult}

\pdfoutput=1

\usepackage{type1cm}        
\usepackage{makeidx}         
\usepackage{graphicx}        
\usepackage{multicol}        
\usepackage{multirow}
\usepackage[bottom]{footmisc}

\usepackage{newtxtext}       %
\usepackage{newtxmath}       

\makeindex             

\usepackage{xcolor}
\usepackage{amsmath,amssymb}
\usepackage{stmaryrd} 
\usepackage{enumitem}
\usepackage{bm,cancel}
\usepackage{accents}
\usepackage{lineno}
\usepackage[T1]{fontenc}

\definecolor{redcol}{rgb}{1.,0.,0.0} 
\definecolor{lnkcol}{rgb}{0.,0.,0.0} 
\usepackage[pdftex,pdfpagelabels=true,bookmarksdepth=3]{hyperref}
\hypersetup{  colorlinks=true,
              citecolor=lnkcol,%
              linkcolor=lnkcol,%
              urlcolor=lnkcol,
            pdfborder=false,
            bookmarks=true,
            bookmarksnumbered=true,
            pdftitle = {On the order reduction of entropy stable DGSEM for the compressible Euler equations},
pdfsubject = {},
pdfauthor = {F.J.Hindenlang,G.J.Gassner},
pdfkeywords = {}
 }
\usepackage[all]{hypcap}  
\usepackage[dvips]{thumbpdf}

\definecolor{chcol}{rgb}{0.4,0.,0.9}

\let\rho\varrho

\renewcommand\vec[1]{\accentset{\,\rightarrow}{#1}}
\newcommand\spacevec[1]{\accentset{\,\rightarrow}{#1}}		
\newcommand\statevec[1]{\mathbf #1}					
\newcommand\acclrvec[1]{\accentset{\,\leftrightarrow}{#1}}	
\newcommand\bigstatevec[1]{\acclrvec{{\mathbf #1}}}		
\newcommand\Ma{\mathrm{Ma}}

\theoremstyle{plain}

\theoremstyle{remark}

\newcommand\tbf[1]{#1}
\newcommand\ordful[1]{\tbf{#1} \textcolor{blue}{$\blacksquare\,$}}
\newcommand\ordred[1]{\tbf{#1} \textcolor{red}{$\boxtimes\,$}}
\newcommand\stent{\rule[-0.1ex]{0pt}{2.5ex}}
\newcommand\cl{\cline{2-11}}

\newcommand\mrtitle[2]{\multirow{12}{*}[-2ex]{\textbf{\begin{tabular}{c} \\ #1 \\ + \\ #2 \end{tabular}}}}
\newcommand\mrm[1]{\multirow{4}{*}[-0.5ex]{#1}}
\newcommand\lrho{$\;\; L_2(\rho)\;\;\;$}
\newcommand\convtable[1]{\begin{scriptsize}\begin{tabular}{|c|c|c|ll|ll|ll|ll|}\hline
 \textbf{Volume disc.} & \multirow{2}{*}{$\;$Mach$\;$} & mesh            & \multicolumn{2}{|c|}{$N=2$} & \multicolumn{2}{|c|}{$N=3$} &\multicolumn{2}{|c|}{$N=4$} & \multicolumn{2}{|c|}{$N=5$} 
 \stent\\ 
 \textbf{+Surface flux} &                               & \hspace{0.5em}level\hspace{0.5em}   & \lrho      & EOC & \lrho     & EOC & \lrho     & EOC& \lrho     & EOC \stent\\\hline\hline
 #1
 \end{tabular}\end{scriptsize}
 }
\newcommand\mrtitleMS[2]{\multirow{4}{*}[-1ex]{\textbf{\begin{tabular}{c} #1 \\ + \\ #2 \end{tabular}}}}
\newcommand\convtableMS[1]{\begin{scriptsize}\begin{tabular}{|c|c|ll|ll|ll|ll|}\hline
 \textbf{Volume disc.} & mesh            & \multicolumn{2}{|c|}{$N=2$} & \multicolumn{2}{|c|}{$N=3$} &\multicolumn{2}{|c|}{$N=4$} & \multicolumn{2}{|c|}{$N=5$} 
 \stent\\ 
 \textbf{+Surface flux} & \hspace{0.5em}level\hspace{0.5em}   & \lrho      & EOC & \lrho     & EOC & \lrho     & EOC& \lrho     & EOC \stent\\\hline\hline
 #1
 \end{tabular}\end{scriptsize}
 }

\begin{document}

\title*{On the order reduction of entropy stable DGSEM for the compressible Euler equations}
\author{Florian J. Hindenlang and Gregor J. Gassner}
\institute{Florian J. Hindenlang \at  Max Planck Institute for Plasma Physics, Boltzmannstra{\ss}e 2, D-85748 Garching, Germany, \email{florian.hindenlang@ipp.mpg.de}
\and Gregor J. Gassner\at Department for Mathematics and Computer Science; Center for Data and Simulation Science, University of Cologne, Cologne, Germany \email{ggassner@math.uni-koeln.de}}
%
%
\maketitle


\abstract{Is the experimental order of convergence lower when using the entropy stable DGSEM-LGL variant? Recently, a debate on the question of the convergence behavior of the entropy stable nodal collocation discontinuous Galerkin spectral element method (DGSEM) with Legendre-Gauss-Lobatto nodes has emerged. Whereas it is well documented that the entropy conservative variant with no additional interface dissipation shows an odd-even behavior when testing its experimental convergence order, the results in the literature are less clear regarding the entropy stable version of the DGSEM-LGL, where explicit Riemann solver type dissipation is added at the element interfaces. We contribute to the ongoing discussion and present numerical experiments for the compressible Euler equations, where we investigate the effect of the choice of the numerical surface flux function. 
In our experiments, it turns out that the choice of the numerical surface flux has an impact on the convergence order. Penalty type numerical fluxes with high dissipation in all waves, such as the LLF and the HLL flux, appear to affect the convergence order negatively for odd polynomial degrees $N$, in contrast to the entropy conserving variant, where even polynomial degrees $N$ are negatively affected. This behavior is more pronounced in low Mach number settings. In contrast, for numerical surface fluxes with less dissipative behavior in the contact wave such as e.g. Roe's flux, the HLLC flux and the entropy conservative flux augmented with 5-wave matrix dissipation, optimal convergence rate of $N+1$ independent of the Mach number is observed.
\keywords{entropy stability, discontinuous Galerkin spectral element method, summation-by-parts, order reduction, compressible Euler equations }
}

\section{Introduction}\label{Sec:Intro}

Discontinuous Galerkin spectral element collocation method (DGSEM) with either Legendre-Gauss or Legendre-Gauss-Lobatto (LGL) nodes (see e.g. \cite{Kopriva:2009nx}) are among the most efficient variants in the class of element based high order methods, such as e.g. discontinuous Galerkin, flux reconstruction, or summation-by-parts (SBP) finite differences. 
In particular, the LGL variant, starting in \cite{gassner_skew_burgers}, turned out to be similar to a SBP finite difference approximation with simultaneous-approximate-term technique (SAT). 
This relationship allowed to construct conservative skew-symmetric approximations, e.g. \cite{gassner_skew_burgers,gassner_kepdg,Winters:2015wj}, and later enabled DGSEM-LGL approximations that are discretely entropy stable, e.g. \cite{carpenter_esdg,PARSANI201588,Niklas-Wintermeyer:2016kq,CHEN2017427,Gassner2018,WINTERMEYER2018447,BOHM2018,Friedrich2018}, and/or kinetic energy preserving \cite{Gassner:2016ye}. These novel variants of nodal split form DG methods feature drastically increased non-linear robustness towards aliasing induced instabilities and favourable properties regarding the simulation of unresolved turbulence, e.g. \cite{FLAD2017782,WINTERS20181}. 

In addition to the very robust dissipative entropy stable versions, it is also possible to construct virtually dissipation free variants by choosing appropriate element interface numerical fluxes. These entropy conserving variants all show an odd-even behavior when experimentally testing the order of convergence, e.g. \cite{gassner_skew_burgers,Winters:2015wj}, where the observed convergence order for even polynomial degrees $N$ is $N$ and for odd $N$ is $N+1$. 
Lately, a discussion emerged in the community, with interesting debates during the recent ICOSAHOM conference in London, where researchers reported non-optimal convergence behavior of the entropy stable DGSEM-LGL even with dissipative numerical surface fluxes, e.g. \cite{CHEN2017427}. 

This paper contributes to this discussion and presents results of an experimental convergence order study for the compressible Euler equations with (i) the standard DGSEM with either Gauss and LGL nodes, (ii) the entropy stable DGSEM with LGL nodes. For these nodal schemes, we test the convergence order with different numerical surface fluxes and report the results depending on the Mach number of the test case.
The remainder of the paper is organized as follows: in the next section we describe the numerical model for our numerical experiments, in Sec.~\ref{Sec:results} we present our observed experimental convergence orders for different configurations and draw our conclusion in  Sec.~\ref{Sec:Concl}.

\section{Numerical Model}\label{Sec:Model}
We consider the compressible Euler equations defined in the domain $\Omega \subset \mathbb{R}^3$ 
\begin{equation}
\label{eq:euler}
{\statevec u_t} + \sum\limits_{i = 1}^3 {\frac{{\partial {\statevec f_{i}}}}{{\partial {x_i}}}} = \statevec 0.
\end{equation}
The state vector contains the conservative variables and the advective flux components are
\begin{equation}\statevec u = \left[ {\begin{array}{*{20}{c}}
  \rho  \\ 
  {\rho \spacevec v} \\ 
  {E} 
\end{array}} \right] = \left[ {\begin{array}{*{20}{c}}
  \rho  \\ 
  {\rho v_1} \\ 
  {\rho v_2} \\ 
  {\rho v_3} \\ 
  {E} 
\end{array}}  \right], \;
\statevec f_{1}  = \left[ {\begin{array}{*{20}c}
   {\rho v_1}  \\
   {\rho v_1^2 + p}  \\
   {\rho v_1\,v_2}  \\
   {\rho v_1\,v_3}  \\
   {( E + p)v_1}  \\

 \end{array} } \right],\; \statevec f_{2}  = \left[ {\begin{array}{*{20}c}
   {\rho v_2}  \\
   {\rho v_2\,v_1}  \\
   {\rho v_2^2  + p}  \\
   {\rho v_2\,v_3}  \\
   {( E + p)v_2}  \\

 \end{array} } \right],\; \statevec f_{3}  = \left[ {\begin{array}{*{20}c}
   {\rho v_3}  \\
   {\rho v_3\,v_1}  \\
   {\rho v_3\,v_2}  \\
   {\rho v_3^2  + p}  \\
   {( E + p) v_3}  \\

 \end{array} } \right].
\end{equation}
Here, $\rho,\,\spacevec{v}=(v_1,v_2,v_3)^T,\,p,\,E$ are the mass density, fluid velocities, pressure and total energy. We close the system with the ideal gas assumption, which relates the total energy and pressure 
\begin{equation}
p = (\gamma-1)\left(E - \frac{1}{2}\rho\left\|\spacevec{v}\right\|^2\right),
\label{eqofstate}
\end{equation}
where $\gamma$ denotes the adiabatic coefficient. 

For our discretization, we subdivide the domain into non-overlapping hexahedral elements. For each element, we define a transfinite mapping to a unit reference space and use this mapping to transform the equations \eqref{eq:euler} from physical to reference space. A weak form is created by taking the inner product of the transformed equation with a test function. We use integration-by-parts for the flux term and approximate the resulting weak form as follows: the conservative variables are approximated by a polynomial in reference space with degree N, interpolated at the Gauss or LGL nodes. The volume fluxes are replaced by a standard interpolation of the non-linear flux function at the same Gauss/LGL nodes (standard DGSEM-Gauss or DGSEM-LGL), see e.g. \cite{Kopriva:2009nx}. For the LGL variant, we are also able to introduce the split form volume integral based on entropy conserving and kinetic energy preserving numerical volume fluxes (Split-DGSEM), e.g. \cite{Gassner:2016ye} and ~\cite{winters2017_ESmatrix}, resulting in either the entropy conserving or entropy stable DGSEM variants, depending on the choice of numerical surface flux. 

\section{Convergence results}\label{Sec:results}
In this section, we compare the convergence of the standard DGSEM and the entropy conservative and entropy stable discretization for different choices of the numerical flux and polynomial degrees $N=2, 3, 4, 5$. 

We choose the test case of a two-dimensional density wave, with a constant pressure and transported with a constant velocity, which was proposed for one-dimensional convergence tests in \cite{CHAN2018346}.  The density evolves as
\begin{equation}
 \rho(x_1,x_2,t)= 1+ 0.1 \sin\left(\pi \left((x_1-v_1 t)+(x_2-v_2 t)\right)\right) \label{eq:density_wave}
\end{equation} 
with a prescribed velocity $(v_1,v_2)$. The pressure is chosen as $p=1/\gamma$ with $\gamma=1.4$, so that the sound speed ranges between $c=0.95\dots 1.05$. 
Thus, by changing the velocity, we change the Mach number of the flow $\Ma=|\spacevec v| / c$. Three Mach numbers are chosen: $\Ma\approx 0.2$ with $(v_1,v_2)=(0.1,0.15)$, $\Ma\approx 1.0$ with $(v_1,v_2)=(0.7,0.65)$ and $\Ma\approx 3.5$ with $(v_1,v_2)=(2.5,2.4)$. The experimental order of convergence (EOC) is computed with the $L_2$ error of the density at $t=1$.

The convergence study is performed with the open source, three-dimensional curvilinear split-form DG framework  FLUXO (\url{www.github.com/project-fluxo}). As the test case is two-dimensional, we use fully periodic cartesian meshes of the domain $[-1,1]^3$ with an equal number of elements in x- and y-directions and always $1$ element in z-direction.
Note that $h_0$ in the convergence tables refers to the coarsest mesh level, which is $4^2$ elements for $N=2,3$ ($h_0=1/2$) and $2^2$ elements for $N=4,5$ ($h_0=1$). 


All simulation results are obtained with an explicit five stage, fourth order accurate low storage Runge-Kutta scheme \cite{Carpenter&Kennedy:1994}, where a stable time step is computed according to the adjustable coefficient CFL$\in(0,1]$, the local maximum wave speed, and the relative grid size, e.g. \cite{gassner2011}. We made sure that the time integrator did not influence the spatial convergence order, by adjusting the CFL number accordingly. 

\begin{table}[htbp!]
\centering 
\convtable{
\mrtitle{DGSEM-Gauss}{HLL} 
 &\mrm{$3.5$} & $h_0/2$  &     1.87e-04  &        3.34  &     8.57e-06  &        4.02  &     1.03e-05  &        5.02  &     6.76e-07  &        6.07    \stent\\
 &            & $h_0/4$  &     2.27e-05  &        3.04  &     5.35e-07  &        4.00  &     3.30e-07  &        4.96  &     1.07e-08  &        5.99    \stent\\
 &            & $h_0/8$  &     2.82e-06  &        3.00  &     3.34e-08  &        4.00  &     1.02e-08  &        5.01  &     1.66e-10  &        6.01    \stent\\
 &            & $h_0/16$ &\tbf{3.53e-07} &\ordful{3.00} &\tbf{2.09e-09} &\ordful{4.00} &\tbf{3.22e-10} &\ordful{4.99} &\tbf{2.60e-12} &\ordful{5.99}\stent\\\cl
 &\mrm{$1.0$} & $h_0/2$  &     2.60e-04  &        2.55  &     5.92e-06  &        4.55  &     1.15e-05  &        4.39  &     6.74e-07  &        6.84    \stent\\ 
 &            & $h_0/4$  &     3.74e-05  &        2.80  &     3.27e-07  &        4.18  &     4.94e-07  &        4.54  &     7.23e-09  &        6.54    \stent\\
 &            & $h_0/8$  &     4.88e-06  &        2.94  &     1.95e-08  &        4.07  &     1.73e-08  &        4.84  &     9.98e-11  &        6.18    \stent\\
 &            & $h_0/16$ &\tbf{6.18e-07} &\ordful{2.98} &\tbf{1.20e-09} &\ordful{4.03} &\tbf{5.60e-10} &\ordful{4.95} &\tbf{1.52e-12} &\ordful{6.04}\stent\\\cl
 &\mrm{$0.2$} & $h_0/2$  &     4.87e-04  &        1.70  &     4.36e-06  &        5.05  &     1.57e-05  &        2.56  &     9.80e-07  &        6.66    \stent\\
 &            & $h_0/4$  &     1.08e-04  &        2.18  &     1.10e-07  &        5.31  &     9.90e-07  &        3.99  &     4.47e-09  &        7.77    \stent\\
 &            & $h_0/8$  &     1.95e-05  &        2.46  &     5.84e-09  &        4.23  &     5.06e-08  &        4.29  &     3.96e-11  &        6.82    \stent\\
 &            & $h_0/16$ &\tbf{3.06e-06} &\ordred{2.67} &\tbf{2.34e-10} &\ordful{4.64} &\tbf{2.26e-09} &\ordred{4.49} &\tbf{4.57e-13} &\ordful{6.44}\stent\\ \hline\hline
\mrtitle{DGSEM-Gauss}{Roe}
 &\mrm{$3.5$} & $h_0/2$  &     1.87e-04  &        3.34  &     8.57e-06  &        4.02  &     1.03e-05  &        5.02  &     6.76e-07  &        6.07   \stent\\
 &            & $h_0/4$  &     2.27e-05  &        3.04  &     5.35e-07  &        4.00  &     3.30e-07  &        4.96  &     1.07e-08  &        5.99    \stent\\
 &            & $h_0/8$  &     2.82e-06  &        3.00  &     3.34e-08  &        4.00  &     1.02e-08  &        5.01  &     1.66e-10  &        6.01    \stent\\
 &            & $h_0/16$ &\tbf{3.53e-07} &\ordful{3.00} &\tbf{2.09e-09} &\ordful{4.00} &\tbf{3.22e-10} &\ordful{4.99} &\tbf{2.60e-12} &\ordful{5.99}   \stent\\\cl
 &\mrm{$1.0$} & $h_0/2$  &     1.82e-04  &        3.07  &     8.76e-06  &        4.04  &     1.11e-05  &        5.07  &     6.95e-07  &        6.08   \stent\\ 
 &            & $h_0/4$  &     2.26e-05  &        3.00  &     5.40e-07  &        4.02  &     3.44e-07  &        5.00  &     1.07e-08  &        6.01    \stent\\
 &            & $h_0/8$  &     2.82e-06  &        3.00  &     3.35e-08  &        4.01  &     1.05e-08  &        5.04  &     1.67e-10  &        6.01    \stent\\
 &            & $h_0/16$ &\tbf{3.53e-07} &\ordful{3.00} &\tbf{2.09e-09} &\ordful{4.00} &\tbf{3.25e-10} &\ordful{5.01} &\tbf{2.61e-12} &\ordful{6.00}\stent\\\cl
 &\mrm{$0.2$} & $h_0/2$  &     2.14e-04  &        2.65  &     1.04e-05  &        3.78  &     1.16e-05  &        4.53  &     7.76e-07  &        5.78   \stent\\
 &            & $h_0/4$  &     2.22e-05  &        3.26  &     5.49e-07  &        4.25  &     3.49e-07  &        5.05  &     1.08e-08  &        6.17    \stent\\
 &            & $h_0/8$  &     2.82e-06  &        2.98  &     3.76e-08  &        3.87  &     1.04e-08  &        5.07  &     1.70e-10  &        5.99    \stent\\
 &            & $h_0/16$ &\tbf{3.53e-07} &\ordful{3.00} &\tbf{2.07e-09} &\ordful{4.19} &\tbf{3.38e-10} &\ordful{4.94} &\tbf{2.64e-12} &\ordful{6.01}   \stent\\ \hline
}
\caption{Experimental order of convergence of $L_2$ error to the exact density \eqref{eq:density_wave}, using the standard DGSEM-Gauss with HLL and Roe fluxes. Full order is marked with \ordful{} ($\gtrsim N+1$)  and an order reduction with \ordred{}.}\label{tab:EOC_density_Gauss_HLL_Roe}
\end{table}

\begin{table}[htbp!]
\centering 
\convtable{
\mrtitle{DGSEM-LGL}{HLL}
 &\mrm{$3.5$} & $h_0/2$  &     1.22e-03  &        3.38  &     3.85e-05  &        4.03  &     4.34e-05  &        4.93  &     2.70e-06  &        5.93    \stent\\
 &            & $h_0/4$  &     1.26e-04  &        3.27  &     2.41e-06  &        4.00  &     1.39e-06  &        4.96  &     4.33e-08  &        5.96    \stent\\
 &            & $h_0/8$  &     1.48e-05  &        3.10  &     1.51e-07  &        4.00  &     4.34e-08  &        5.01  &     6.67e-10  &        6.02    \stent\\
 &            & $h_0/16$ &\tbf{1.81e-06} &\ordful{3.03} &\tbf{9.42e-09} &\ordful{4.00} &\tbf{1.36e-09} &\ordful{4.99} &\tbf{1.05e-11} &\ordful{5.98} \stent\\\cl
 &\mrm{$1.0$} & $h_0/2$  &     1.04e-03  &        2.30  &     3.44e-05  &        4.44  &     4.11e-05  &        4.89  &     2.89e-06  &        6.33    \stent\\ 
 &            & $h_0/4$  &     1.58e-04  &        2.72  &     1.88e-06  &        4.20  &     1.77e-06  &        4.54  &     3.85e-08  &        6.23    \stent\\
 &            & $h_0/8$  &     2.11e-05  &        2.91  &     1.16e-07  &        4.02  &     6.25e-08  &        4.82  &     5.37e-10  &        6.16    \stent\\
 &            & $h_0/16$ &\tbf{2.69e-06} &\ordful{2.97} &\tbf{7.18e-09} &\ordful{4.01} &\tbf{2.05e-09} &\ordful{4.93} &\tbf{8.36e-12} &\ordful{6.00} \stent\\\cl
 &\mrm{$0.2$} & $h_0/2$  &     1.20e-03  &        1.96  &     4.86e-05  &        3.51  &     5.16e-05  &        3.91  &     3.85e-06  &        5.31    \stent\\
 &            & $h_0/4$  &     2.72e-04  &        2.14  &     1.85e-06  &        4.71  &     2.99e-06  &        4.11  &     4.25e-08  &        6.50    \stent\\
 &            & $h_0/8$  &     5.57e-05  &        2.29  &     1.21e-07  &        3.94  &     1.56e-07  &        4.26  &     5.45e-10  &        6.28    \stent\\
 &            & $h_0/16$ &\tbf{1.01e-05} &\ordred{2.47} &\tbf{5.95e-09} &\ordful{4.34} &\tbf{7.19e-09} &\ordred{4.44} &\tbf{7.11e-12} &\ordful{6.26}   \stent\\ \hline\hline
\mrtitle{DGSEM-LGL}{Roe}
 &\mrm{$3.5$} & $h_0/2$  &     1.22e-03  &        3.38  &     3.85e-05  &        4.03  &     4.34e-05  &        4.93  &     2.70e-06  &        5.93   \stent\\
 &            & $h_0/4$  &     1.26e-04  &        3.27  &     2.41e-06  &        4.00  &     1.39e-06  &        4.96  &     4.33e-08  &        5.96    \stent\\
 &            & $h_0/8$  &     1.48e-05  &        3.10  &     1.51e-07  &        4.00  &     4.34e-08  &        5.01  &     6.67e-10  &        6.02    \stent\\
 &            & $h_0/16$ &\tbf{1.81e-06} &\ordful{3.03} &\tbf{9.42e-09} &\ordful{4.00} &\tbf{1.36e-09} &\ordful{4.99} &\tbf{1.05e-11} &\ordful{5.98} \stent\\\cl
 &\mrm{$1.0$} & $h_0/2$  &     9.17e-04  &        2.86  &     3.96e-05  &        3.94  &     4.41e-05  &        4.94  &     2.76e-06  &        6.02   \stent\\ 
 &            & $h_0/4$  &     1.15e-04  &        2.99  &     2.41e-06  &        4.04  &     1.47e-06  &        4.90  &     4.44e-08  &        5.96    \stent\\
 &            & $h_0/8$  &     1.44e-05  &        3.00  &     1.51e-07  &        4.00  &     4.38e-08  &        5.07  &     6.84e-10  &        6.02    \stent\\
 &            & $h_0/16$ &\tbf{1.80e-06} &\ordful{3.00} &\tbf{9.42e-09} &\ordful{4.00} &\tbf{1.35e-09} &\ordful{5.02} &\tbf{1.08e-11} &\ordful{5.99} \stent\\\cl
 &\mrm{$0.2$} & $h_0/2$  &     9.26e-04  &        2.35  &     4.63e-05  &        3.45  &     4.26e-05  &        4.27  &     2.97e-06  &        5.37   \stent\\
 &            & $h_0/4$  &     1.19e-04  &        2.96  &     2.40e-06  &        4.27  &     1.59e-06  &        4.74  &     4.48e-08  &        6.05    \stent\\
 &            & $h_0/8$  &     1.43e-05  &        3.06  &     1.57e-07  &        3.93  &     4.34e-08  &        5.19  &     6.83e-10  &        6.04    \stent\\
 &            & $h_0/16$ &\tbf{1.80e-06} &\ordful{2.99} &\tbf{9.37e-09} &\ordful{4.07} &\tbf{1.49e-09} &\ordful{4.86} &\tbf{1.09e-11} &\ordful{5.97} \stent\\\hline
}
\caption{Experimental order of convergence of $L_2$ error to the exact density \eqref{eq:density_wave}, using DGSEM-GL with HLL and Roe fluxes. Full order is marked with \ordful{} ($\gtrsim N+1$)  and an order reduction with \ordred{}.}\label{tab:EOC_density_GL_HLL_Roe}
\end{table}

\subsection{Standard DGSEM}

The convergence of the standard DGSEM with Gauss-Legendre nodes (DGSEM-Gauss) and with Legendre-Gauss-Lobatto (DGSEM-LGL) is shown in Table~\ref{tab:EOC_density_Gauss_HLL_Roe} and Table~\ref{tab:EOC_density_GL_HLL_Roe}, for the three Mach numbers and two choices of the numerical flux, namely the HLL and the Roe flux. The results of the local Lax-Friedrichs (LLF) and the HLLC fluxes are reported in the Appendix, as HLL is similar to LLF, and HLLC is the same as Roe, see Table~\ref{tab:EOC_density_Gauss_LLF_HLLC} and Table~\ref{tab:EOC_density_GL_LLF_HLLC}. 

For the HLL flux and the low Mach number $\Ma=0.2$, we observe an odd-even behavior with an \emph{order reduction} for even polynomial degrees $N=2,4$. Also for $\Ma=1.0$, the convergence for even degrees is slightly affected, whereas for the high Mach number, all fluxes converge with full order. Comparing the $L_2$ errors of the finest mesh for HLL and Roe for the low Mach number, HLL is   less accurate for $N=2,4$ and   more accurate for $N=3,5$.

All numerical fluxes are approximate Riemann solvers, but the LLF and HLL only use the maximum wave speeds, whereas the HLLC and Roe also take the contact wave into account, and therefore keep the full order of the scheme for all Mach numbers for this test case. 

\begin{table}[h!]
\centering
\convtable{
\mrtitle{Split-DGSEM}{ECKEP}
 &\mrm{$3.5$} & $h_0/2$  &     1.62e-03  &        4.23  &     8.45e-05  &        2.77  &     5.80e-05  &        5.88  &     4.63e-06  &        5.63    \stent\\
 &            & $h_0/4$  &     1.30e-04  &        3.64  &     7.14e-06  &        3.56  &     1.60e-06  &        5.18  &     8.38e-08  &        5.79    \stent\\
 &            & $h_0/8$  &     1.05e-05  &\ordful{3.62} &     7.90e-07  &        3.18  &     4.56e-08  &        5.13  &     3.59e-09  &        4.54    \stent\\
 &            & $h_0/16$ &\tbf{1.69e-06} &\ordred{2.64} &\tbf{9.58e-08} &\ordred{3.04} &\tbf{1.23e-09} &\ordful{5.21} &\tbf{1.18e-10} &\ordred{4.93} \stent\\\cl
 &\mrm{$1.0$} & $h_0/2$  &     1.41e-03  &        3.89  &     9.45e-05  &        2.34  &     7.71e-05  &        5.86  &     3.49e-06  &        6.00    \stent\\ 
 &            & $h_0/4$  &     1.25e-04  &        3.49  &     1.26e-05  &        2.90  &     1.98e-06  &        5.28  &     5.00e-08  &        6.12    \stent\\
 &            & $h_0/8$  &     1.48e-05  &        3.09  &     1.60e-06  &        2.98  &     4.06e-08  &        5.61  &     1.22e-09  &        5.36    \stent\\
 &            & $h_0/16$ &\tbf{1.32e-06} &\ordful{3.48} &\tbf{2.01e-07} &\ordred{3.00} &\tbf{1.13e-09} &\ordful{5.17} &\tbf{3.85e-11} &\ordred{4.98}   \stent\\\cl
 &\mrm{$0.2$} & $h_0/2$  &     1.13e-03  &        2.26  &     8.03e-05  &        2.98  &     5.95e-05  &        4.14  &     4.11e-06  &        5.14    \stent\\
 &            & $h_0/4$  &     1.13e-04  &        3.32  &     1.02e-05  &        2.97  &     1.85e-06  &        5.01  &     1.33e-07  &        4.95    \stent\\
 &            & $h_0/8$  &     1.12e-05  &\ordful{3.34} &     1.29e-06  &        2.99  &     4.02e-08  &        5.52  &     4.21e-09  &        4.98    \stent\\
 &            & $h_0/16$ &\tbf{1.98e-06} &\ordred{2.50} &\tbf{1.61e-07} &\ordred{3.00} &\tbf{1.43e-09} &\ordful{4.81} &\tbf{1.32e-10} &\ordred{5.00}   \stent\\ \hline\hline
\mrtitle{Split-DGSEM}{HLL}
 &\mrm{$3.5$} & $h_0/2$  &     1.23e-03  &        3.36  &     3.88e-05  &        4.06  &     4.49e-05  &        4.94  &     3.01e-06  &        5.97    \stent\\
 &            & $h_0/4$  &     1.27e-04  &        3.27  &     2.42e-06  &        4.00  &     1.43e-06  &        4.97  &     4.84e-08  &        5.96    \stent\\
 &            & $h_0/8$  &     1.48e-05  &        3.10  &     1.51e-07  &        4.00  &     4.44e-08  &        5.01  &     7.48e-10  &        6.02    \stent\\
 &            & $h_0/16$ &\tbf{1.81e-06} &\ordful{3.03} &\tbf{9.46e-09} &\ordful{4.00} &\tbf{1.40e-09} &\ordful{4.99} &\tbf{1.18e-11} &\ordful{5.98} \stent\\\cl
 &\mrm{$1.0$} & $h_0/2$  &     1.04e-03  &        2.30  &     3.48e-05  &        4.44  &     4.16e-05  &        4.90  &     3.48e-06  &        6.10    \stent\\ 
 &            & $h_0/4$  &     1.58e-04  &        2.72  &     1.88e-06  &        4.21  &     1.79e-06  &        4.54  &     4.66e-08  &        6.23    \stent\\
 &            & $h_0/8$  &     2.11e-05  &        2.90  &     1.16e-07  &        4.02  &     6.39e-08  &        4.81  &     6.08e-10  &        6.26    \stent\\
 &            & $h_0/16$ &\tbf{2.69e-06} &\ordful{2.97} &\tbf{7.20e-09} &\ordful{4.01} &\tbf{2.11e-09} &\ordful{4.92} &\tbf{9.28e-12} &\ordful{6.03}   \stent\\\cl
 &\mrm{$0.2$} & $h_0/2$  &     1.20e-03  &        1.96  &     4.90e-05  &        3.51  &     5.23e-05  &        3.92  &     4.47e-06  &        5.16    \stent\\
 &            & $h_0/4$  &     2.72e-04  &        2.14  &     1.87e-06  &        4.71  &     3.02e-06  &        4.12  &     6.14e-08  &        6.19    \stent\\
 &            & $h_0/8$  &     5.58e-05  &        2.29  &     1.21e-07  &        3.95  &     1.57e-07  &        4.26  &     6.18e-10  &        6.64    \stent\\
 &            & $h_0/16$ &\tbf{1.01e-05} &\ordred{2.47} &\tbf{5.98e-09} &\ordful{4.34} &\tbf{7.28e-09} &\ordred{4.43} &\tbf{8.20e-12} &\ordful{6.23}   \stent\\ \hline\hline
\mrtitle{Split-DGSEM}{ECKEP-Roe}
 &\mrm{$3.5$} & $h_0/2$  &     1.23e-03  &        3.36  &     3.88e-05  &        4.06  &     4.49e-05  &        4.94  &     3.01e-06  &        5.97   \stent\\
 &            & $h_0/4$  &     1.27e-04  &        3.27  &     2.42e-06  &        4.00  &     1.43e-06  &        4.97  &     4.84e-08  &        5.96    \stent\\
 &            & $h_0/8$  &     1.48e-05  &        3.10  &     1.51e-07  &        4.00  &     4.44e-08  &        5.01  &     7.48e-10  &        6.02    \stent\\
 &            & $h_0/16$ &\tbf{1.81e-06} &\ordful{3.03} &\tbf{9.46e-09} &\ordful{4.00} &\tbf{1.40e-09} &\ordful{4.99} &\tbf{1.18e-11} &\ordful{5.98}   \stent\\\cl
 &\mrm{$1.0$} & $h_0/2$  &     9.18e-04  &        2.86  &     3.98e-05  &        3.94  &     4.51e-05  &        4.92  &     3.08e-06  &        5.95   \stent\\ 
 &            & $h_0/4$  &     1.15e-04  &        2.99  &     2.42e-06  &        4.04  &     1.51e-06  &        4.90  &     4.97e-08  &        5.95    \stent\\
 &            & $h_0/8$  &     1.44e-05  &        3.00  &     1.51e-07  &        4.00  &     4.49e-08  &        5.07  &     7.66e-10  &        6.02    \stent\\
 &            & $h_0/16$ &\tbf{1.80e-06} &\ordful{3.00} &\tbf{9.46e-09} &\ordful{4.00} &\tbf{1.38e-09} &\ordful{5.02} &\tbf{1.21e-11} &\ordful{5.99}   \stent\\\cl
 &\mrm{$0.2$} & $h_0/2$  &     9.26e-04  &        2.35  &     4.65e-05  &        3.45  &     4.37e-05  &        4.26  &     3.26e-06  &        5.34   \stent\\
 &            & $h_0/4$  &     1.19e-04  &        2.96  &     2.40e-06  &        4.27  &     1.63e-06  &        4.74  &     5.03e-08  &        6.02    \stent\\
 &            & $h_0/8$  &     1.43e-05  &        3.06  &     1.58e-07  &        3.93  &     4.45e-08  &        5.19  &     7.66e-10  &        6.04    \stent\\
 &            & $h_0/16$ &\tbf{1.80e-06} &\ordful{2.99} &\tbf{9.40e-09} &\ordful{4.07} &\tbf{1.53e-09} &\ordful{4.86} &\tbf{1.22e-11} &\ordful{5.97}   \stent\\\hline
}
\caption{Experimental order of convergence of $L_2$ error to the exact density \eqref{eq:density_wave}, using entropy conservative ECKEP flux and entropy stable HLL and ECKEP-Roe fluxes. Full order is marked with \ordful{} ($\gtrsim N+1$)  and an order reduction with \ordred{}. }\label{tab:EOC_density_wave_ECKEP_ECKEP_HLL_Roe}
\end{table}

\subsection{Entropy conservative and entropy stable DGSEM}

Now, we investigate the order reduction of the entropy conservative and entropy stable discretizations. Here, the standard DGSEM volume integral is replaced by split-form formulation (Split-DGSEM) using a two-point entropy conservative and kinetic energy preserving flux ECKEP. If we choose the ECKEP flux at the surface, we get an entropy-conserving scheme. For entropy stability, we can use the LLF or HLL flux directly at the surface, or use the ECKEP flux and add a dissipation term, which must still satisfy the entropy inequality condition. In Winters et al.~\cite{winters2017_ESmatrix}, such dissipation terms are carefully derived, using either only the maximum wave speed (LLF-type) or incorporating all waves (Roe-type), which we will refer to as ECKEP-LLF and ECKEP-Roe fluxes. 

In Table~\ref{tab:EOC_density_wave_ECKEP_ECKEP_HLL_Roe}, we summarize the convergence of the dissipation-free ECKEP flux, the HLL and ECKEP-Roe flux. The results for LLF and ECKEP-LF fluxes are found in the Appendix in Table~\ref{tab:EOC_density_wave_ECKEP_LF_ECKEP-LF}, as they have the same convergence and error levels as the HLL flux. 
As expected, the dissipation-free surface flux (ECKEP) produces an order reduction for all Mach numbers for $N=3,5$, and for $N=2$ full order not kept in the last refinement step. 

If we simply use the HLL flux, we have an entropy stable scheme, but an order reduction for $N=2,4$ can be observed for the low Mach number flow, analogously to the standard DGSEM-LGL scheme. Interestingly, the odd-even behavior switches between entropy conserving and entropy stable fluxes. 

The ECKEP-Roe entropy stable flux accounts for all waves of the Riemann problem and adjusts the dissipation for each wave accordingly, which gives full order convergence for all Mach numbers.


%

\section{Conclusions}\label{Sec:Concl}
In this work, we report the convergence of standard DGSEM Gauss and Gauss-Lobatto schemes to entropy conservative (EC) and entropy stable (ES) DGSEM schemes for the Euler equations, as there have been findings of order reduction for EC and ES schemes.  We choose a simple density transport test case on a periodic domain and investigate the influence of the Mach number of the transport velocity. 

The EC scheme is dissipation free and an order reduction is observed by the convergence study presented here, confirming many similar observations found in literature. 
We also confirm that the ES scheme can have an order reduction for low Mach numbers, but only if the entropy stable numerical flux relies on simple approximate Riemann solvers such as local Lax-Friedrichs or HLL. If all waves are accounted for in the dissipation term of the entropy stable flux as presented in \cite{winters2017_ESmatrix}, the full order is observed for all Mach numbers. IN addition, we reproduce the same behavior for the standard DGSEM Gauss and Gauss-Lobatto schemes, where the LLF and HLL fluxes suffer from order reduction at low Mach number, and HLLC and Roe fluxes have full order for all Mach numbers. 

We want to emphasize that the present convergence study should be seen merely as an observation, confirming that the numerical flux can have strong influence on the convergence order for both the standard DGSEM and the entropy stable DGSEM. Also, we stress that in our tests the order reduction is related to the form of the dissipation term in the numerical surface flux and is not related to the insufficient integration precision of the LGL-quadrature. 

In the authors' experience, a convergence study using a manufactured solution technique can be misleading, as full convergence order is found independent of the choice of numerical flux. Hence, the introduction of a source term to balance the prescribed solution overcomes possible deficiencies of the surface fluxes, showing the limit of the manufactured solution technique in this context. In the Appendix, the convergence results of a manufactured solution are reported.

\section*{Acknowledgements}
Gregor Gassner thanks the European Research Council for funding through the ERC Starting Grant ``An Exascale aware and Un-crashable Space-Time-Adaptive Discontinuous Spectral Element Solver for Non-Linear Conservation Laws'' (\texttt{Extreme}), ERC grant agreement no. 714487.
Florian Hindenlang thanks Eric Sonnendr\"ucker and the Max-Planck Institute for Plasma Physics in Garching for their constant support. 
We would also like all participants of the ICOSAHOM 2018 for the valuable discussions on the topic of entropy stable schemes, which motivated this work.

\bibliography{dakBib}

\section*{Appendix}

\subsection*{Additional convergence results}
In this section, we present additional convergence results of the density wave tst case for the DGSEM-Gauss and DGSEM-LGL with LLF and HLLC fluxes in Table~\ref{tab:EOC_density_Gauss_LLF_HLLC} and Table~\ref{tab:EOC_density_GL_LLF_HLLC}, and also the entropy stable schemes with LLF and ECKEP-LLF fluxes in Table~\ref{tab:EOC_density_wave_ECKEP_LF_ECKEP-LF}. The results for LLF-type fluxes behave like the HLL flux, and for the HLLC flux like the Roe-type fluxes presented in Table~\ref{tab:EOC_density_wave_ECKEP_ECKEP_HLL_Roe}.

\begin{table}[htbp!]
\centering 
\convtable{
\mrtitle{DGSEM-Gauss}{LLF}
 &\mrm{$3.5$} & $h_0/2$  &     2.42e-04  &        2.97  &     6.43e-06  &        4.41  &     1.05e-05  &        4.51  &     6.68e-07  &        6.60    \stent\\
 &            & $h_0/4$  &     3.24e-05  &        2.90  &     3.71e-07  &        4.11  &     4.32e-07  &        4.60  &     8.04e-09  &        6.38    \stent\\
 &            & $h_0/8$  &     4.15e-06  &        2.96  &     2.27e-08  &        4.03  &     1.47e-08  &        4.87  &     1.15e-10  &        6.13    \stent\\
 &            & $h_0/16$ &\tbf{5.22e-07} &\ordful{2.99} &\tbf{1.41e-09} &\ordful{4.01} &\tbf{4.73e-10} &\ordful{4.96} &\tbf{1.77e-12} &\ordful{6.02} \stent\\\cl 
 &\mrm{$1.0$} & $h_0/2$  &     3.13e-04  &        2.29  &     4.59e-06  &        4.84  &     1.18e-05  &        3.90  &     6.69e-07  &        7.40    \stent\\ 
 &            & $h_0/4$  &     5.30e-05  &        2.56  &     2.25e-07  &        4.35  &     6.08e-07  &        4.28  &     5.56e-09  &        6.91    \stent\\
 &            & $h_0/8$  &     7.43e-06  &        2.83  &     1.29e-08  &        4.12  &     2.47e-08  &        4.62  &     6.79e-11  &        6.35    \stent\\
 &            & $h_0/16$ &\tbf{9.61e-07} &\ordful{2.95} &\tbf{7.65e-10} &\ordful{4.07} &\tbf{8.53e-10} &\ordful{4.85} &\tbf{9.97e-13} &\ordful{6.09}   \stent\\\cl 
 &\mrm{$0.2$} & $h_0/2$  &     4.95e-04  &        1.69  &     4.33e-06  &        5.07  &     1.58e-05  &        2.53  &     9.88e-07  &        6.68    \stent\\
 &            & $h_0/4$  &     1.12e-04  &        2.15  &     1.06e-07  &        5.35  &     1.01e-06  &        3.97  &     4.46e-09  &        7.79    \stent\\
 &            & $h_0/8$  &     2.06e-05  &        2.44  &     5.47e-09  &        4.28  &     5.23e-08  &        4.27  &     3.86e-11  &        6.85    \stent\\
 &            & $h_0/16$ &\tbf{3.29e-06} &\ordred{2.65} &\tbf{2.15e-10} &\ordful{4.67} &\tbf{2.38e-09} &\ordred{4.46} &\tbf{4.35e-13} &\ordful{6.47}   \stent\\ \hline\hline
\mrtitle{DGSEM-Gauss}{HLLC} 
 &\mrm{$3.5$} & $h_0/2$  &     1.87e-04  &        3.34  &     8.57e-06  &        4.02  &     1.03e-05  &        5.02  &     6.76e-07  &        6.07   \stent\\
 &            & $h_0/4$  &     2.27e-05  &        3.04  &     5.35e-07  &        4.00  &     3.30e-07  &        4.96  &     1.07e-08  &        5.99    \stent\\
 &            & $h_0/8$  &     2.82e-06  &        3.00  &     3.34e-08  &        4.00  &     1.02e-08  &        5.01  &     1.66e-10  &        6.01    \stent\\
 &            & $h_0/16$ &\tbf{3.53e-07} &\ordful{3.00} &\tbf{2.09e-09} &\ordful{4.00} &\tbf{3.22e-10} &\ordful{4.99} &\tbf{2.60e-12} &\ordful{5.99}   \stent\\\cl
 &\mrm{$1.0$} & $h_0/2$  &     1.82e-04  &        3.07  &     8.76e-06  &        4.04  &     1.11e-05  &        5.07  &     6.95e-07  &        6.08   \stent\\ 
 &            & $h_0/4$  &     2.26e-05  &        3.00  &     5.40e-07  &        4.02  &     3.44e-07  &        5.00  &     1.07e-08  &        6.01    \stent\\
 &            & $h_0/8$  &     2.82e-06  &        3.00  &     3.35e-08  &        4.01  &     1.05e-08  &        5.04  &     1.67e-10  &        6.01    \stent\\
 &            & $h_0/16$ &\tbf{3.53e-07} &\ordful{3.00} &\tbf{2.09e-09} &\ordful{4.00} &\tbf{3.25e-10} &\ordful{5.01} &\tbf{2.61e-12} &\ordful{6.00}   \stent\\\cl
 &\mrm{$0.2$} & $h_0/2$  &     2.14e-04  &        2.65  &     1.04e-05  &        3.78  &     1.16e-05  &        4.53  &     7.76e-07  &        5.78   \stent\\
 &            & $h_0/4$  &     2.22e-05  &        3.26  &     5.49e-07  &        4.25  &     3.49e-07  &        5.05  &     1.08e-08  &        6.17    \stent\\
 &            & $h_0/8$  &     2.82e-06  &        2.98  &     3.76e-08  &        3.87  &     1.04e-08  &        5.07  &     1.70e-10  &        5.99    \stent\\
 &            & $h_0/16$ &\tbf{3.53e-07} &\ordful{3.00} &\tbf{2.07e-09} &\ordful{4.19} &\tbf{3.38e-10} &\ordful{4.94} &\tbf{2.64e-12} &\ordful{6.01}   \stent\\\hline
}
\caption{Experimental order of convergence of $L_2$ error to the exact density \eqref{eq:density_wave}, using DGSEM-Gauss with LLF and HLLC fluxes. Full order is marked with \ordful{} ($\gtrsim N+1$)  and an order reduction with \ordred{}.}\label{tab:EOC_density_Gauss_LLF_HLLC}
\end{table}

\begin{table}[htbp!]
\centering 
\convtable{
\mrtitle{DGSEM-LGL}{LLF}
 &\mrm{$3.5$} & $h_0/2$  &     1.23e-03  &        3.11  &     3.44e-05  &        4.36  &     4.02e-05  &        4.82  &     2.91e-06  &        6.14    \stent\\
 &            & $h_0/4$  &     1.51e-04  &        3.02  &     2.00e-06  &        4.10  &     1.61e-06  &        4.64  &     3.80e-08  &        6.26    \stent\\
 &            & $h_0/8$  &     1.89e-05  &        3.00  &     1.23e-07  &        4.02  &     5.55e-08  &        4.86  &     5.46e-10  &        6.12    \stent\\
 &            & $h_0/16$ &\tbf{2.36e-06} &\ordful{3.00} &\tbf{7.66e-09} &\ordful{4.01} &\tbf{1.79e-09} &\ordful{4.95} &\tbf{8.54e-12} &\ordful{6.00} \stent\\\cl
 &\mrm{$1.0$} & $h_0/2$  &     1.10e-03  &        1.97  &     3.15e-05  &        4.83  &     4.02e-05  &        4.93  &     2.85e-06  &        6.73    \stent\\ 
 &            & $h_0/4$  &     2.03e-04  &        2.43  &     1.68e-06  &        4.23  &     1.98e-06  &        4.34  &     3.68e-08  &        6.28    \stent\\
 &            & $h_0/8$  &     2.99e-05  &        2.76  &     1.03e-07  &        4.02  &     8.35e-08  &        4.57  &     4.86e-10  &        6.24    \stent\\
 &            & $h_0/16$ &\tbf{3.94e-06} &\ordful{2.93} &\tbf{6.39e-09} &\ordful{4.01} &\tbf{2.96e-09} &\ordful{4.82} &\tbf{7.53e-12} &\ordful{6.01}   \stent\\\cl
 &\mrm{$0.2$} & $h_0/2$  &     1.21e-03  &        1.95  &     4.88e-05  &        3.51  &     5.20e-05  &        3.90  &     3.89e-06  &        5.31    \stent\\
 &            & $h_0/4$  &     2.77e-04  &        2.12  &     1.86e-06  &        4.72  &     3.05e-06  &        4.09  &     4.27e-08  &        6.51    \stent\\
 &            & $h_0/8$  &     5.78e-05  &        2.26  &     1.21e-07  &        3.94  &     1.61e-07  &        4.24  &     5.47e-10  &        6.29    \stent\\
 &            & $h_0/16$ &\tbf{1.06e-05} &\ordred{2.44} &\tbf{5.98e-09} &\ordful{4.34} &\tbf{7.57e-09} &\ordred{4.41} &\tbf{7.16e-12} &\ordful{6.26}   \stent\\ \hline\hline
\mrtitle{DGSEM-LGL}{HLLC}
 &\mrm{$3.5$} & $h_0/2$  &     1.22e-03  &        3.38  &     3.85e-05  &        4.03  &     4.34e-05  &        4.93  &     2.70e-06  &        5.93   \stent\\
 &            & $h_0/4$  &     1.26e-04  &        3.27  &     2.41e-06  &        4.00  &     1.39e-06  &        4.96  &     4.33e-08  &        5.96    \stent\\
 &            & $h_0/8$  &     1.48e-05  &        3.10  &     1.51e-07  &        4.00  &     4.34e-08  &        5.01  &     6.67e-10  &        6.02    \stent\\
 &            & $h_0/16$ &\tbf{1.81e-06} &\ordful{3.03} &\tbf{9.42e-09} &\ordful{4.00} &\tbf{1.36e-09} &\ordful{4.99} &\tbf{1.05e-11} &\ordful{5.98}   \stent\\\cl
 &\mrm{$1.0$} & $h_0/2$  &     9.17e-04  &        2.86  &     3.96e-05  &        3.94  &     4.41e-05  &        4.94  &     2.76e-06  &        6.02   \stent\\ 
 &            & $h_0/4$  &     1.15e-04  &        2.99  &     2.41e-06  &        4.04  &     1.47e-06  &        4.90  &     4.44e-08  &        5.96    \stent\\
 &            & $h_0/8$  &     1.44e-05  &        3.00  &     1.51e-07  &        4.00  &     4.38e-08  &        5.07  &     6.84e-10  &        6.02    \stent\\
 &            & $h_0/16$ &\tbf{1.80e-06} &\ordful{3.00} &\tbf{9.42e-09} &\ordful{4.00} &\tbf{1.35e-09} &\ordful{5.02} &\tbf{1.08e-11} &\ordful{5.99}   \stent\\\cl
 &\mrm{$0.2$} & $h_0/2$  &     9.26e-04  &        2.35  &     4.63e-05  &        3.45  &     4.26e-05  &        4.27  &     2.97e-06  &        5.37   \stent\\
 &            & $h_0/4$  &     1.19e-04  &        2.96  &     2.40e-06  &        4.27  &     1.59e-06  &        4.74  &     4.48e-08  &        6.05    \stent\\
 &            & $h_0/8$  &     1.43e-05  &        3.06  &     1.57e-07  &        3.93  &     4.34e-08  &        5.19  &     6.83e-10  &        6.04    \stent\\
 &            & $h_0/16$ &\tbf{1.80e-06} &\ordful{2.99} &\tbf{9.37e-09} &\ordful{4.07} &\tbf{1.49e-09} &\ordful{4.86} &\tbf{1.09e-11} &\ordful{5.97}   \stent\\\hline
}
\caption{Experimental order of convergence of $L_2$ error to the exact density \eqref{eq:density_wave}, using DGSEM-GL with LLF and HLLC fluxes. Full order is marked with \ordful{} ($\gtrsim N+1$)  and an order reduction with \ordred{}.}\label{tab:EOC_density_GL_LLF_HLLC}
\end{table}

\begin{table}[htbp!]
\centering
\convtable{
\mrtitle{Split-DGSEM}{LLF}
 &\mrm{$3.5$} & $h_0/2$  &     1.24e-03  &        3.11  &     3.49e-05  &        4.37  &     4.16e-05  &        4.85  &     3.37e-06  &        6.05    \stent\\
 &            & $h_0/4$  &     1.52e-04  &        3.03  &     2.01e-06  &        4.12  &     1.64e-06  &        4.67  &     4.49e-08  &        6.23    \stent\\
 &            & $h_0/8$  &     1.89e-05  &        3.00  &     1.23e-07  &        4.03  &     5.67e-08  &        4.85  &     6.21e-10  &        6.18    \stent\\
 &            & $h_0/16$ &\tbf{2.36e-06} &\ordful{3.00} &\tbf{7.68e-09} &\ordful{4.01} &\tbf{1.84e-09} &\ordful{4.95} &\tbf{9.58e-12} &\ordful{6.02} \stent\\\cl
 &\mrm{$1.0$} & $h_0/2$  &     1.09e-03  &        1.97  &     3.21e-05  &        4.82  &     4.07e-05  &        4.94  &     3.77e-06  &        6.34    \stent\\ 
 &            & $h_0/4$  &     2.03e-04  &        2.43  &     1.69e-06  &        4.25  &     1.99e-06  &        4.35  &     4.69e-08  &        6.33    \stent\\
 &            & $h_0/8$  &     2.99e-05  &        2.76  &     1.04e-07  &        4.03  &     8.46e-08  &        4.56  &     5.57e-10  &        6.39    \stent\\
 &            & $h_0/16$ &\tbf{3.94e-06} &\ordful{2.93} &\tbf{6.41e-09} &\ordful{4.02} &\tbf{3.03e-09} &\ordful{4.80} &\tbf{8.38e-12} &\ordful{6.05}   \stent\\\cl
 &\mrm{$0.2$} & $h_0/2$  &     1.21e-03  &        1.95  &     4.92e-05  &        3.50  &     5.26e-05  &        3.91  &     4.52e-06  &        5.15    \stent\\
 &            & $h_0/4$  &     2.77e-04  &        2.12  &     1.88e-06  &        4.71  &     3.07e-06  &        4.10  &     6.21e-08  &        6.19    \stent\\
 &            & $h_0/8$  &     5.79e-05  &        2.26  &     1.22e-07  &        3.95  &     1.63e-07  &        4.24  &     6.21e-10  &        6.64    \stent\\
 &            & $h_0/16$ &\tbf{1.06e-05} &\ordred{2.44} &\tbf{6.01e-09} &\ordful{4.34} &\tbf{7.66e-09} &\ordred{4.41} &\tbf{8.26e-12} &\ordful{6.23}   \stent\\ \hline\hline
\mrtitle{Split-DGSEM}{ECKEP-LLF}
 &\mrm{$3.5$} & $h_0/2$  &     1.24e-03  &        3.10  &     3.49e-05  &        4.37  &     4.16e-05  &        4.84  &     3.38e-06  &        6.06    \stent\\
 &            & $h_0/4$  &     1.52e-04  &        3.03  &     2.01e-06  &        4.12  &     1.64e-06  &        4.67  &     4.49e-08  &        6.23    \stent\\
 &            & $h_0/8$  &     1.89e-05  &        3.01  &     1.23e-07  &        4.03  &     5.68e-08  &        4.85  &     6.21e-10  &        6.18    \stent\\
 &            & $h_0/16$ &\tbf{2.36e-06} &\ordful{3.00} &\tbf{7.68e-09} &\ordful{4.01} &\tbf{1.84e-09} &\ordful{4.95} &\tbf{9.58e-12} &\ordful{6.02} \stent\\\cl
 &\mrm{$1.0$} & $h_0/2$  &     1.10e-03  &        1.97  &     3.20e-05  &        4.83  &     4.08e-05  &        4.94  &     3.77e-06  &        6.36    \stent\\ 
 &            & $h_0/4$  &     2.04e-04  &        2.43  &     1.69e-06  &        4.25  &     2.00e-06  &        4.35  &     4.69e-08  &        6.33    \stent\\
 &            & $h_0/8$  &     3.00e-05  &        2.76  &     1.04e-07  &        4.03  &     8.49e-08  &        4.56  &     5.57e-10  &        6.40    \stent\\
 &            & $h_0/16$ &\tbf{3.95e-06} &\ordful{2.93} &\tbf{6.41e-09} &\ordful{4.02} &\tbf{3.04e-09} &\ordful{4.81} &\tbf{8.38e-12} &\ordful{6.05}   \stent\\\cl
 &\mrm{$0.2$} & $h_0/2$  &     1.21e-03  &        1.95  &     4.93e-05  &        3.51  &     5.27e-05  &        3.91  &     4.53e-06  &        5.15    \stent\\
 &            & $h_0/4$  &     2.78e-04  &        2.12  &     1.88e-06  &        4.71  &     3.08e-06  &        4.10  &     6.21e-08  &        6.19    \stent\\
 &            & $h_0/8$  &     5.79e-05  &        2.26  &     1.22e-07  &        3.95  &     1.63e-07  &        4.24  &     6.21e-10  &        6.65    \stent\\
 &            & $h_0/16$ &\tbf{1.07e-05} &\ordred{2.44} &\tbf{6.01e-09} &\ordful{4.34} &\tbf{7.67e-09} &\ordred{4.41} &\tbf{8.26e-12} &\ordful{6.23}  \stent\\\hline 
}
\caption{Experimental order of convergence of $L_2$ error to the exact density \eqref{eq:density_wave}, using entropy stable LLF and ECKEP-LLF flux. Full order is marked with \ordful{} ($\gtrsim N+1$)  and an order reduction with \ordred{}. }\label{tab:EOC_density_wave_ECKEP_LF_ECKEP-LF}
\end{table}

\subsection*{Manufactured solution with source term}\label{Sec:Conv_MS}

Here, we run a convergence test with the method of manufactured solutions. To do so, we assume a two-dimensional solution of the form
\begin{equation}
\begin{split}
\statevec{u} = &\left[ \rho \,, \rho v_1\,, \rho v_2 \,, \rho v_3 \,, E  \right]^T=\left[  g\,,g\,,g\,,0\,,g^2  \right]^T \\
&\text{ with } \; g = g(x_1,x_2,t) = 0.5 \sin(2\pi(x_1+x_2-t))+2.
\label{eq:MS} 
\end{split}
\end{equation}
Note that the average Mach number in the domain is $\Ma=0.8$. Inserting \eqref{eq:MS} into the Euler equations, and using the fact that spatial and time derivatives are $g'=\partial_{x_1} g=\partial_{x_2}g=-\partial_t g$, we get an additional residual
\begin{equation}
\statevec{u}_t + \vec{\nabla} \cdot \bigstatevec{f}(\statevec{u}) = \begin{pmatrix} g' \\ (3\gamma-2)g' + 2(\gamma-1) g g' \\ (3\gamma-2)g' + 2(\gamma-1) g g' \\ 0 \\ (6\gamma-2) g' + 2(2\gamma-1) g g'  \end{pmatrix}
\label{eq:MS_residual}
\end{equation}
To solve the inhomogeneous problem, we subtract the residual from the approximate solution in each Runge-Kutta step. Moreover, we run the test case up to the final time $t\!=\!1.0$. 

In the convergence results for the standard DGSEM Gauss and Gauss-Lobatto, we see that the LLF flux still leads to an order reduction for $N=2,4$, whereas full order is found for the HLL, HLLC and Roe fluxes, see Table~\ref{tab:EOC_mansol_Gauss} and Table~\ref{tab:EOC_mansol_LGL}.  

In Table~\ref{tab:EOC_mansol_ECKEP} the entropy conservative scheme shows again an order reduction for $N=3,5$, and the LLF-Type dissipation too, for $N=2,4$, and for this test case, all entropy stable schemes exhibit full order.

\begin{table}[htbp!]
\centering 
\convtableMS{
\mrtitleMS{DGSEM-Gauss}{LLF}  & $h_0/2$  &     2.30e-03  &        2.20  &     4.54e-05  &        5.34  &     1.13e-04  &        6.35  &     4.52e-05  &        4.38    \stent\\
                              & $h_0/4$  &     4.81e-04  &        2.25  &     1.99e-06  &        4.52  &     4.78e-06  &        4.56  &     2.37e-07  &        7.58    \stent\\
                              & $h_0/8$  &     9.48e-05  &        2.34  &     1.02e-07  &        4.28  &     2.88e-07  &        4.05  &     1.53e-09  &        7.28    \stent\\
                              & $h_0/16$ &\tbf{1.57e-05} &\ordred{2.60} &\tbf{6.25e-09} &\ordful{4.03} &\tbf{1.73e-08} &\ordred{4.05} &\tbf{1.57e-11} &\ordful{6.61}   \stent\\ \hline\hline
\mrtitleMS{DGSEM-Gauss}{HLL}  & $h_0/2$  &     1.24e-03  &        2.84  &     5.46e-05  &        4.35  &     1.32e-04  &        5.22  &     1.47e-05  &        6.30   \stent\\
                              & $h_0/4$  &     1.17e-04  &        3.41  &     3.36e-06  &        4.02  &     2.89e-06  &        5.51  &     1.44e-07  &        6.67    \stent\\
                              & $h_0/8$  &     1.41e-05  &        3.04  &     1.85e-07  &        4.18  &     7.14e-08  &        5.34  &     1.74e-09  &        6.37    \stent\\
                              & $h_0/16$ &\tbf{1.76e-06} &\ordful{3.00} &\tbf{1.07e-08} &\ordful{4.11} &\tbf{2.15e-09} &\ordful{5.05} &\tbf{2.24e-11} &\ordful{6.28}   \stent\\\hline
\mrtitleMS{DGSEM-Gauss}{HLLC} & $h_0/2$  &     1.24e-03  &        2.84  &     5.46e-05  &        4.35  &     1.32e-04  &        5.22  &     1.47e-05  &        6.30   \stent\\
                              & $h_0/4$  &     1.17e-04  &        3.41  &     3.36e-06  &        4.02  &     2.89e-06  &        5.51  &     1.44e-07  &        6.67    \stent\\
                              & $h_0/8$  &     1.41e-05  &        3.04  &     1.85e-07  &        4.18  &     7.14e-08  &        5.34  &     1.74e-09  &        6.37    \stent\\
                              & $h_0/16$ &\tbf{1.76e-06} &\ordful{3.00} &\tbf{1.07e-08} &\ordful{4.11} &\tbf{2.15e-09} &\ordful{5.05} &\tbf{2.24e-11} &\ordful{6.28}   \stent\\\hline
\mrtitleMS{DGSEM-Gauss}{Roe}  & $h_0/2$  &     1.24e-03  &        2.84  &     5.46e-05  &        4.35  &     1.32e-04  &        5.22  &     1.47e-05  &        6.30   \stent\\
                              & $h_0/4$  &     1.17e-04  &        3.41  &     3.36e-06  &        4.02  &     2.89e-06  &        5.51  &     1.44e-07  &        6.67    \stent\\
                              & $h_0/8$  &     1.41e-05  &        3.04  &     1.85e-07  &        4.18  &     7.14e-08  &        5.34  &     1.74e-09  &        6.37    \stent\\
                              & $h_0/16$ &\tbf{1.76e-06} &\ordful{3.00} &\tbf{1.07e-08} &\ordful{4.11} &\tbf{2.15e-09} &\ordful{5.05} &\tbf{2.24e-11} &\ordful{6.28}   \stent\\\hline
}
\caption{Experimental order of convergence of $L_2$ error of density for the manufactured solution \eqref{eq:MS}, using DGSEM-Gauss with LLF, HLL, HLLC and Roe fluxes. Full order is marked with \ordful{} ($\gtrsim N+1$)  and an order reduction with \ordred{}.}\label{tab:EOC_mansol_Gauss}
\end{table}

\begin{table}[htbp!]
\centering 
\convtableMS{
\mrtitleMS{DGSEM-LGL}{LLF}    & $h_0/2$  &     7.36e-03  &        2.85  &     3.15e-04  &        4.38  &     5.69e-04  &        5.78  &     9.33e-05  &        5.34    \stent\\
                              & $h_0/4$  &     1.33e-03  &        2.47  &     1.43e-05  &        4.46  &     2.04e-05  &        4.80  &     9.27e-07  &        6.65    \stent\\
                              & $h_0/8$  &     2.79e-04  &        2.25  &     7.99e-07  &        4.16  &     8.81e-07  &        4.54  &     9.34e-09  &        6.63    \stent\\
                              & $h_0/16$ &\tbf{5.31e-05} &\ordred{2.39} &\tbf{4.72e-08} &\ordful{4.08} &\tbf{5.94e-08} &\ordred{3.89} &\tbf{1.37e-10} &\ordful{6.09}   \stent\\ \hline\hline
\mrtitleMS{DGSEM-LGL}{HLL}    & $h_0/2$  &     5.32e-03  &        3.16  &     2.52e-04  &        4.00  &     3.84e-04  &        5.36  &     4.33e-05  &        6.24   \stent\\
                              & $h_0/4$  &     5.99e-04  &        3.15  &     1.38e-05  &        4.19  &     1.43e-05  &        4.75  &     4.58e-07  &        6.56    \stent\\
                              & $h_0/8$  &     7.25e-05  &        3.05  &     7.69e-07  &        4.17  &     2.92e-07  &        5.61  &     7.08e-09  &        6.02    \stent\\
                              & $h_0/16$ &\tbf{9.02e-06} &\ordful{3.01} &\tbf{4.74e-08} &\ordful{4.02} &\tbf{7.77e-09} &\ordful{5.23} &\tbf{1.10e-10} &\ordful{6.01}   \stent\\\hline
\mrtitleMS{DGSEM-LGL}{HLLC}   & $h_0/2$  &     5.32e-03  &        3.16  &     2.52e-04  &        4.00  &     3.84e-04  &        5.36  &     4.33e-05  &        6.24   \stent\\
                              & $h_0/4$  &     5.99e-04  &        3.15  &     1.38e-05  &        4.19  &     1.43e-05  &        4.75  &     4.58e-07  &        6.56    \stent\\
                              & $h_0/8$  &     7.25e-05  &        3.05  &     7.69e-07  &        4.17  &     2.92e-07  &        5.61  &     7.08e-09  &        6.02    \stent\\
                              & $h_0/16$ &\tbf{9.02e-06} &\ordful{3.01} &\tbf{4.74e-08} &\ordful{4.02} &\tbf{7.77e-09} &\ordful{5.23} &\tbf{1.10e-10} &\ordful{6.01}   \stent\\\hline
\mrtitleMS{DGSEM-LGL}{Roe}    & $h_0/2$  &     5.32e-03  &        3.16  &     2.52e-04  &        4.00  &     3.84e-04  &        5.36  &     4.33e-05  &        6.24   \stent\\
                              & $h_0/4$  &     5.99e-04  &        3.15  &     1.38e-05  &        4.19  &     1.43e-05  &        4.75  &     4.58e-07  &        6.56    \stent\\
                              & $h_0/8$  &     7.25e-05  &        3.05  &     7.69e-07  &        4.17  &     2.92e-07  &        5.61  &     7.08e-09  &        6.02    \stent\\
                              & $h_0/16$ &\tbf{9.02e-06} &\ordful{3.01} &\tbf{4.74e-08} &\ordful{4.02} &\tbf{7.77e-09} &\ordful{5.23} &\tbf{1.10e-10} &\ordful{6.01}   \stent\\\hline
}
\caption{Experimental order of convergence of $L_2$ error of density for the manufactured solution \eqref{eq:MS}, using DGSEM-LGL with LLF, HLL, HLLC and Roe fluxes. Full order is marked with \ordful{} ($\gtrsim N+1$)  and an order reduction with \ordred{}.}\label{tab:EOC_mansol_LGL}
\end{table}

\begin{table}[htbp!]
\centering 
\convtableMS{
\mrtitleMS{Split-DGSEM}{ECKEP}      & $h_0/2$  &     1.31e-02  &        3.53  &     1.28e-03  &        2.88  &     5.62e-03  &        2.49  &     1.12e-03  &        4.05    \stent\\
                              & $h_0/4$  &     1.30e-03  &        3.34  &     1.13e-04  &        3.49  &     3.12e-04  &        4.17  &     5.97e-06  &        7.56    \stent\\
                              & $h_0/8$  &     1.23e-04  &        3.40  &     1.24e-05  &        3.20  &     2.00e-06  &        7.29  &     9.45e-08  &        5.98    \stent\\
                              & $h_0/16$ &\tbf{1.76e-05} &\ordful{2.80} &\tbf{1.67e-06} &\ordred{2.89} &\tbf{3.41e-08} &\ordful{5.87} &\tbf{3.17e-09} &\ordred{4.90}   \stent\\ \hline\hline
\mrtitleMS{Split-DGSEM}{LLF}        & $h_0/2$  &     7.60e-03  &        3.13  &     3.70e-04  &        3.91  &     6.56e-04  &        4.44  &     1.28e-04  &        4.75    \stent\\
                              & $h_0/4$  &     1.63e-03  &        2.22  &     1.90e-05  &        4.29  &     3.09e-05  &        4.41  &     2.01e-06  &        6.00    \stent\\
                              & $h_0/8$  &     3.40e-04  &        2.26  &     9.89e-07  &        4.26  &     1.83e-06  &        4.08  &     1.80e-08  &        6.80    \stent\\
                              & $h_0/16$ &\tbf{6.14e-05} &\ordred{2.47} &\tbf{6.41e-08} &\ordful{3.95} &\tbf{9.15e-08} &\ordred{4.32} &\tbf{2.17e-10} &\ordful{6.37}   \stent\\ \hline\hline
\mrtitleMS{Split-DGSEM}{ECKEP-LLF}  & $h_0/2$  &     7.66e-03  &        3.21  &     3.74e-04  &        3.98  &     6.79e-04  &        4.45  &     1.33e-04  &        4.75   \stent\\
                              & $h_0/4$  &     1.64e-03  &        2.23  &     1.90e-05  &        4.30  &     3.08e-05  &        4.46  &     2.06e-06  &        6.01    \stent\\
                              & $h_0/8$  &     3.41e-04  &        2.26  &     9.88e-07  &        4.26  &     1.83e-06  &        4.07  &     1.79e-08  &        6.85    \stent\\
                              & $h_0/16$ &\tbf{6.16e-05} &\ordred{2.47} &\tbf{6.40e-08} &\ordful{3.95} &\tbf{9.17e-08} &\ordred{4.32} &\tbf{2.17e-10} &\ordful{6.37}   \stent\\\hline
\mrtitleMS{Split-DGSEM}{HLL}        & $h_0/2$  &     5.82e-03  &        3.10  &     3.01e-04  &        3.85  &     5.11e-04  &        4.35  &     7.06e-05  &        5.27   \stent\\
                              & $h_0/4$  &     7.06e-04  &        3.04  &     2.04e-05  &        3.88  &     1.67e-05  &        4.94  &     1.08e-06  &        6.03    \stent\\
                              & $h_0/8$  &     8.63e-05  &        3.03  &     1.16e-06  &        4.14  &     5.08e-07  &        5.04  &     1.67e-08  &        6.02    \stent\\
                              & $h_0/16$ &\tbf{1.08e-05} &\ordful{3.00} &\tbf{7.20e-08} &\ordful{4.01} &\tbf{1.62e-08} &\ordful{4.97} &\tbf{2.64e-10} &\ordful{5.98}   \stent\\\hline
\mrtitleMS{Split-DGSEM}{ECKEP-Roe}  & $h_0/2$  &     5.81e-03  &        3.11  &     3.01e-04  &        3.85  &     5.11e-04  &        4.35  &     7.06e-05  &        5.27   \stent\\
                              & $h_0/4$  &     7.06e-04  &        3.04  &     2.04e-05  &        3.88  &     1.67e-05  &        4.94  &     1.08e-06  &        6.03    \stent\\
                              & $h_0/8$  &     8.63e-05  &        3.03  &     1.16e-06  &        4.14  &     5.08e-07  &        5.04  &     1.67e-08  &        6.02    \stent\\
                              & $h_0/16$ &\tbf{1.08e-05} &\ordful{3.00} &\tbf{7.20e-08} &\ordful{4.01} &\tbf{1.62e-08} &\ordful{4.97} &\tbf{2.64e-10} &\ordful{5.98}   \stent\\\hline
}
\caption{Experimental order of convergence of $L_2$ error of density for the manufactured solution \eqref{eq:MS}, using entropy conservative and entropy stable schemes. Full order is marked with \ordful{} ($\gtrsim N+1$)  and an order reduction with \ordred{}.}\label{tab:EOC_mansol_ECKEP}
\end{table}

\end{document}